\documentclass{amsart}

\usepackage{amsmath,xypic,leftidx}
\usepackage{xspace}
\usepackage[psamsfonts]{amssymb}
\usepackage[latin1]{inputenc}
\usepackage{graphicx,color}

\usepackage{hyperref}

\theoremstyle{remark}

\theoremstyle{plain}

\newtheorem{proposition}{\bf Proposition}

\newtheorem{theoreme}{\bf Théorème}
\newtheorem{lemme}{\bf Lemme}

\def\cal{\mathcal}

\def\C{{\mathbb C}}
\def\R{{\mathbb R}}

\def\P{{{\mathbb P}}}
\def\T{{\cal T}}
\def\Z{{\cal Z}}

\def\loc{{\rm loc}}
\def\ddc{{\rm dd}^{\rm c}}

\def\eps{{\varepsilon}}

\title{Courants dynamiques pluripolaires}
\begin{author}[X.~Buff]{Xavier Buff}
\email{xavier.buff$@$math.univ-toulouse.fr}
\address{ %
  Institut de Mathématiques de Toulouse\\
  Universit\'e Paul Sabatier\\
  118, route de Narbonne \\
  31062 Toulouse Cedex \\
  France }
\end{author}

\begin{document}

\begin{abstract}
On montre l'existence d'applications rationnelles $f:\P^k\to \P^k$ telles que
\begin{itemize}
\item $f$ est algébriquement stable~: pour tout $n\geq 0$, $\deg f^{\circ n} =(\deg f)^n$, 

\item il existe un unique courant positif fermé $T$ de bidegré $(1,1)$ vérifiant $f^* T=d\cdot T$ et $\int_{\P^k} T\wedge \omega^{k-1}=1$ où $\omega$ est la forme de Fubini-Study sur $\P^k$ et 

\item $T$ est pluripolaire~: il existe un ensemble pluripolaire $X\subset \P^k$ tel que $\int_{X}T\wedge \omega^{k-1}=1$. 
\end{itemize}
\end{abstract}

\maketitle

On se donne un entier $k\geq 2$. On note $\pi:\C^{k+1}-\{0\}\to \P^k$ la projection canonique. Soit $\omega$ la forme de Fubini-Study sur $\P^k$ définie par $\pi^*\omega = \ddc \log \|\cdot\|$.

Rappelons qu'une fonction $u:\C^{k+1}\to [-\infty,+\infty)$, non identiquement $-\infty$, est psh (plurisousharmonique) si $u$ est semi-continue supérieurement et si la restriction de $u$ à toute droite complexe est soit sousharmonique, soit identiquement $-\infty$. Alors, $u\in L^1_\loc(\C^{k+1})$ et le courant  $\ddc u$ est positif. 

Nous dirons qu'une fonction  $g:\P^k\to[-\infty,+\infty)$ est $\omega$-psh lorsque la fonction $u:\C^{k+1}\to [-\infty,+\infty)$ définie par $u(0)=-\infty$ et $u(z)=\log \|z\|+g\circ \pi(z)$ si $z\neq 0$ est une fonction psh sur $\C^{k+1}$. En particulier, 
\begin{itemize}
\item $g$ est semi-continue supérieurement, 
\item $g$ n'est pas identiquement $-\infty$, 
\item $g\in L^1(\P^k)$,
\item $\omega+\ddc g$ est un courant positif et 
\item $\ddc u = \pi^*(\omega+\ddc g)$.
\end{itemize}

On note $\T$ l'ensemble des courants positifs fermés et de bidegré $(1,1)$ sur $\P^k$, équipé de la topologie faible. Si $T\in \T$, la distribution $T\wedge \omega^{k-1}$ s'identifie, via le théorème de représentation de Riesz, avec une mesure borélienne sur $\P^k$. La masse de $T$ est 
\[\|T\|=\int_{\P^k} T\wedge \omega^{k-1}.\] Plus généralement, si $X\subseteq \P^k$ est un ensemble mesurable, on note
\[\|T\|_X = \int_{X} T\wedge \omega^{k-1}.\] 
L'ensemble des courants $T\in \T$ de masse $\|T\|=1$ est compact. 
Lorsque $\|T\|=1$, il existe une fonction $\omega$-psh $g:\P^k\to[-\infty,+\infty)$ telle que
$\omega+\ddc g = T$. On dit que $g$ est un potentiel associé à $T$. Notons que si $g_1$ et $g_2$ sont deux potentiels associés au même courant $T$, alors $g_1-g_2$ est une fonction pluriharmonique sur $\P^k$, donc constante.

Soit $f:\P^k\to \P^k$ une application rationnelle de degré algébrique $d\geq 2$ et soit $F:\C^{k+1}\to \C^{k+1}$ une application vérifiant $f\circ  \pi=\pi\circ F$ et dont les coordonnées sont des polynômes homogènes de degré $d$. Dans tout l'article, les applications rationnelles $f$ considérées seront supposées {\em dominantes}, c'est-à-dire que le Jacobien de $F$ est non identiquement nul. 
L'application $f:\P^k\to \P^k$ induit alors une application continue $f^*:\T\to \T$ et l'on a $\|f^* T\|=d\cdot \|T\|$. 
Il est facile de voir qu'il existe un courant $T\in \T$ tel que
\[f^*T = d\cdot T\quad \text{et}\quad \|T\|=1.\]
On peut par exemple considérer la suite de courants  $(T_N)$ définis par
\[S_0=\omega,\quad S_{n+1}=\frac{1}{d}f^* S_n\quad \text{et}\quad T_N= \frac{1}{N}\sum_{n=0}^{N-1} S_n.\]
Par récurrence sur $n\geq 0$, on a $\|S_n\|=1$ et donc $\|T_N\|=1$ pour tout $N\geq 1$. 
Par passage à la limite sur
\[\frac{1}{d}f^* T_N = T_N + \frac{1}{N}(S_N-S_0),\]
on obtient que toute valeur d'adhérence $T$ de la suite $T_N$ satisfait   
\[f^*T = d\cdot T\quad \text{et}\quad \|T\|=1.\]
Dans le cas général, on ne sait pas s'il y a unicité d'un tel courant $T$. 

Rappelons qu'un ensemble $X\subset \P^k$ est {\em pluripolaire} s'il existe une fonction $\omega$-psh $g$ sur $\P^k$ telle que $X=\{g=-\infty\}$. 
Nous dirons qu'un courant $T\in \T$ est {\em pluripolaire} s'il existe un ensemble pluripolaire $X\subset \P^k$ tel que $\|T\|_X  = \|T\|$
(toute la masse de $T$ se concentre sur un ensemble pluripolaire).

Une application rationnelle $f:\P^k\to \P^k$ est {\em algébriquement stable} si le degré algébrique de $f$ et le degré algébrique de $f^{\circ n}$ sont reliés par $\deg f^{\circ n} =(\deg f)^n$ pour tout $n\geq 0$. 
Notre but est de montrer l'existence d'applications rationnelles algébriquement stables $f:\P^k\to \P^k$   pour lesquelles il existe un unique courant $T\in \T$ vérifiant $f^*T = d\cdot T$ et $\|T\|=1$,  et pour lesquelles ce courant $T$ est pluripolaire.

\begin{theoreme}
Soient $\bigl(f_t:\P^k\to \P^k)_{t\in [0,1]}$ une famille d'applications rationnelles de degré $d$ et $\Z\subset [0,1]$ un ensemble dénombrable dense tels que
\begin{itemize}
\item il existe une famille $(F_t:\C^{k+1}\to \C^{k+1})_{t\in [0,1]}$ vérifiant $f_t\circ \pi = \pi \circ F_t$ et dépendant continuement de $t\in[0,1]$ pour la topologie de la convergence localement uniforme,

\item l'application $F_t$ fixe un point $z_t\in \C^{k+1}-\{0\}$ qui dépend continuement de $t$ et

\item  l'application $f_t$ est algébriquement stable si et seulement si $t\in [0,1]-\Z$. 
\end{itemize}
Alors, il existe un ensemble $A\subset [0,1]-\Z$ gras au sens de Baire tel que pour tout $t\in A$, 
\begin{itemize}
\item il existe un unique courant $T\in \T$, vérifiant $f_t^* T=d\cdot T$ et $\|T\|=1$ et

\item $T$ est pluripolaire. 
\end{itemize}
\end{theoreme}

Un exemple d'une telle famille d'applications rationnelles $\bigl(f_t:\P^k\to \P^k)_{t\in [0,1]}$ avec $k=2$ et $d=2$ est donné par Charles Favre \cite{f}. Dans son exemple, les applications $f_t$ sont birationnelles. 
En appliquant le théorème à $(f_t)$ et $(f_t^{-1})$ on montre que l'ensemble des $t\in [0,1]-\Z$ pour lesquels les courants de Green $T^+_{f_t}$ et $T^-_{f_t}$ sont pluripolaires est gras au sens de Baire. 

La variante suivante est due à Jeffrey Diller et Vincent Guedj \cite{dg}. Considérons la famille d'applications birationnelles $f_t:\P^2\to \P^2$ définies par
\[F_t(x,y,z) = \bigl(bcx(-cx+acy+z),acy(x-ay+abz),abz(bcx+y-bz)\bigr)\]
avec
\[a=i, \quad b=-2e^{i\pi/4}e^{i\pi t}\quad \text{et}\quad c=\frac{1}{2}e^{i\pi/4}e^{i\pi t}.\]
L'inverse de $f_t$ est une application de la même forme où $a$, $b$, $c$, sont remplacés par $a^{-1}$, $b^{-1}$ et $c^{-1}$.  
Les points d'indétermination de $f_t$ sont $[a:1:0]$, $[0:b:1]$ et $[1:0:c]$ et ceux de $f_t^{-1}$ sont 
$[a^{-1}:1:0]$, $[0:b^{-1}:1]$ et $[1:0:c^{-1}]$. 
Les applications $f_t$ et $f_t^{-1}$ fixent globalement chacune des droites $\{x=0\}$, $\{y=0\}$ et $\{z=0\}$. 

Dans la droite $\{x=0\}$, le point d'indétermination de $f_t$ est $[0:b:1]$ et le point d'indétermination de $f_t^{-1}$ est $[0:1/b:1]$. La restriction de $f_t$ à cette droite est $[0:y:1]\mapsto [0:iy/4:1]$. \'Etant donné que $|1/b|<|b|$, l'orbite de $[0:1/b:1]$ sous itération de $f_t$ ne rencontre jamais $[0:b:1]$ et  l'orbite de $[0:b:1]$ sous itération de $f_t^{-1}$ ne rencontre jamais $[0:1/b:1]$. 

Dans la droite $\{y=0\}$, le point d'indétermination de $f_t$ est $[1:0:c]$ et le point d'indétermination de $f_t^{-1}$ est $[1:0:1/c]$. La restriction de $f_t$ à cette droite est $[1:0:z]\mapsto [1:0:-4iz]$. \'Etant donné que $|1/c|>|c|$, l'orbite de $[1:0:1/c]$ sous itération de $f_t$ ne rencontre jamais $[1:0:c]$ et l'orbite de $[1:0:c]$ sous itération de $f_t^{-1}$ ne rencontre jamais $[1:0:1/c]$. 

Dans la droite $\{z=0\}$, le point d'indétermination de $f_t$ est $[i:1:0]$ et le point d'indétermination de $f_t^{-1}$ est $[-i:1:0]$. La restriction de $f_t$ à cette droite est $[x:1:0]\mapsto [e^{i2\pi t} x:1:0]$, c'est-à-dire une rotation d'angle $t$ tours. Les points $[i:1:0]$ et $[-i:1:0]$ sont deux points opposés d'un cercle invariant par $f_t$. 
\begin{itemize}
\item 
Si $t$ est irrationnel, l'orbite de $[-i:1:0]$ sous itération de $f_t$ ne rencontre donc pas le point $[i:1:0]$ et l'orbite de $[i:0:1]$ sous itération de $f_t^{-1}$ ne rencontre pas $[-i:0:1]$. Dans ce cas les applications $f_t$ et $f_t^{-1}$ sont algébriquement stables. 
\item
En revanche, si $t=(2p+1)/(2q)$ est rationnel à dénominateur pair, alors le $q$-ième itéré de $f_t$ est la rotation d'angle $1/2$ tour qui envoie $[-i:1:0]$ sur $[i:1:0]$. De même le $q$-ième itéré de $f_t^{-1}$ est la rotation d'angle $1/2$ tour qui envoie $[i:1:0]$ sur $[-i:1:0]$. Les applications $f_t$ et $f_t^{-1}$ ne sont pas algébriquement stables. 
\end{itemize}

Enfin, $F_t(0,0,z)=(0,0,-ab^2z^2)$ et donc, l'application $F_t$ fixe le point 
\[z_t=\left(0,0,-\frac{1}{ab^2}\right) = \left(0,0,e^{i2\pi t}/4\right)\in \C^3-\{0\}\]
qui dépend continuement de $t$. 

Jeffrey Diller et Vincent Guedj \cite{dg} montrent que dans le cas de leur exemple, lorsque $f$ est algébriquement stable, il est toujours possible de définir une mesure produit
\[\mu_f=T^+_f\wedge T^-_f.\]
D'après Jeffrey Diller, Romain Dujardin et Vincent Guedj \cite{ddg}, lorsque les courants $T^+_f$ et $T^-_f$ sont pluripolaires, l'application $f$ n'est pas d'{\em énergie finie} et la mesure $\mu_f$ charge totalement un ensemble  pluripolaire. Notre résultat répond donc à une question posée dans \cite{ddg} concernant l'existence d'applications rationnelles de petit degré topologique et d'énergie infinie. Bien que dans ce cas, les techniques exposées dans \cite{ddg} ne s'appliquent pas, il reste néanmoins intéressant d'étudier les propriétés dynamiques de la mesure $\mu_f$.

\section{Courants de Green dynamiques\label{sec_green}}

La construction suivante est due à Nessim Sibony \cite{s}. Soit $f:\P^k\to \P^k$ une application rationnelle de degré $d\geq 2$ et soit $F:\C^{k+1}\to \C^{k+1}$ une application polynômiale homogène relevant $f$, c'est-à-dire telle que $f\circ \pi = \pi\circ F$. Sans perte de généralité, quitte à multiplier $F$ par une constante $C$, on peut supposer que $\|F(z)\|< \|z\|^d$ quand $\|z\|<1$.

On définit la fonction $\omega$-psh $g^{(0)}$ par
\[g^{(0)}\bigl( \pi(z)\bigr)= \frac{1}{d} \log \|F(z)\| - \log \|z\|.\]
On vérifie facilement que si la fonction $\omega$-psh $h$ est un potentiel pour un courant $T\in \T$ de masse $1$, alors $g^{(0)}+ d^{-1}\cdot h\circ f$ est un potentiel pour $d^{-1}f^*T$. 
Trouver un courant $T\in \T$ tel que $f^*T = d\cdot T$ et $\|T\|=1$ revient donc à trouver une fonction $\omega$-psh $h$ telle que $g^{(0)}+d^{-1}\cdot h\circ f =h+c$ avec $c\in \R$. 
Quitte à remplacer $h$ par $h-dc/(d-1)$, cela se ramène à determiner les solutions $\omega$-psh $h$ de l'équation
\begin{equation}\label{eq_cohomologique} h-\frac{1}{d} h\circ f = g^{(0)}.\end{equation}

Nous avons vu qu'il existe toujours un courant $T\in \T$ tel que $f^*T = d\cdot T$ et $\|T\|=1$, donc une solution $\omega$-psh $h$ de l'équation (\ref{eq_cohomologique}) avec $T=\omega+\ddc h$. 
On peut également produire une solution de la manière suivante. On considère la suite de fonctions $\omega$-psh $(g^{(n)})_{n\geq 0}$ définie récursivement pour $n\geq 1$ par
\[g^{(n)} = g^{(0)}+\frac{1}{d}g^{(n-1)}\circ f\]
de sorte que $T^{(n)}=\omega+\ddc g^{(n)}$ satisfasse
\[T^{(0)}=\omega \quad \text{et}\quad T^{(n)} = \frac{1}{d}f^* T^{(n-1)}.\]
Alors,
\[g^{(n)} = \sum_{m=0}^{n-1} \frac{1}{d^m} g^{(0)} \circ f^{\circ m} = h-\frac{1}{d^{n+1}} h\circ f^{\circ n+1}.\]
\'Etant donné que $g^{(0)}$ est négative, la suite $(g^{(n)})_{n\geq 0}$ est décroissante. La fonction $h$ est majorée sur $\P^k$. La suite  $(g^{(n)})_{n\geq 0}$ est donc  minorée par $h$. Elle converge simplement et dans $L^1(\P^k)$ vers une fonction $\omega$-psh $g$. \'Etant donné que 
\[g^{(n)}-\frac{1}{d}g^{(n-1)}\circ f = g^{(0)},\] cette limite $g$ est solution de l'équation (\ref{eq_cohomologique}). 
On note $T_f$ le courant de Green dynamique défini par
\[T_f= \omega+\ddc g.\]

D'après ce que nous venons de voir, toute solution $h$ de l'équation (\ref{eq_cohomologique}) vérifie 
$h\leq g+c$ avec $c\in \R$. 

Notons de plus que si l'on fait l'hypothèse que $F$ fixe un point $z_0\in \C^{k+1}-\{0\}$ (c'est le cas dans le cadre du théorème) et si l'on note $x_0=\pi(z_0)$, alors  
\[g^{(0)}(x_0)=\frac{1}{d}\log \|z_0\| -\log \|z_0\|\]
et on a le contrôle uniforme suivant: 
\[g(x_0) = \sum_{n=0}^{+\infty}\frac{1}{d^n} g^{(0)}\circ f^{\circ n}(x_0) =  \sum_{n=0}^{+\infty}\frac{1}{d^n} g^{(0)}(x_0)=\frac{d}{d-1}g^{(0)}(x_0)=-\log \|z_0\|.\] 

\section{Le cas non algébriquement stable}

Dans cette partie, on suppose que $f:\P^k\to \P^k$ est une application rationnelle de degré $d$ qui n'est pas algébriquement stable. On suppose que $F:\C^{k+1}\to \C^{k+1}$ est une application dont les coordonnées sont des polynômes homogènes de degré $d$ telle que $f\circ \pi = \pi \circ F$ et  $\|F(z)\|<\|z\|^d$ dès que $\|z\|<1$. 
Comme dans la partie précédente, on définit une suite décroissante de fonctions $\omega$-psh $(g^{(n)})_{n\geq 0}$  qui convergent simplement et dans $L^1(\P^k)$ vers une fonction $\omega$-psh $g$ et on pose $T_f=\omega+\ddc g$. 

\begin{proposition}[Sibony \cite{s}]\label{prop_sibony}
Si $f$ est algébriquement instable, le courant de Green dynamique $T_f$ est porté par une réunion dénombrable d'ensembles algébriques. En particulier, il est pluripolaire.
\end{proposition}

\begin{proof}[Preuve] 
\'Etant donné que $f$ n'est pas algébriquement stable,  si $n$ est suffisamment grand, $\deg(f^{\circ n})<d^n$. Cela signifie que les coordonnées de $F^{\circ n}$ ont un facteur commun. On note $H^{(n)}$ le plus grand diviseur commun des coordonnées de $F^{\circ n}$ de sorte que
\[F^{\circ n} = H^{(n)}\cdot \check F^{(n)}\quad
\text{avec}\quad \deg f^{\circ n} = \deg \check F^{(n)}.\]
Alors, 
\[g^{(n)} \bigl(\pi(z)\bigr)= \frac{1}{d^n}\log \|H^{(n)}(z)\| + \frac{1}{d^n}\log \|\check F^{(n)}(z)\|-\log \|z\|.\]
En particulier, le courant $\omega+\ddc g^{(n)}$ charge l'ensemble algébrique $X_n$ défini par l'équation $H_n=0$ avec 
\[\|\omega+\ddc g^{(n)}\|_{X_n}=\frac{1}{d^n}\deg(H^{(n)})=1-\frac{\deg(f^{\circ n})}{d^n}.\] 
De plus, pour $m\geq 0$,
\[F^{\circ n+m} = (H^{(n)})^{d^m}\cdot F^{\circ m}\circ  \check F^{(n)},\]
ce qui montre que 
$(H^{(n)})^{d^m}$ divise $H^{(n+m)}$. 
On en déduit que pour tout $m\geq 0$, le courant $\omega+\ddc g^{(n+m)}$ charge l'ensemble algébrique $X_n$ avec 
\[\|\omega+\ddc g^{(n+m)}\|_{X_n}\geq 1-\frac{\deg(f^{\circ n})}{d^n}.\] 
Quand $m\to +\infty$, la suite de courants $\omega+\ddc g^{(n+m)}$ converge faiblement vers le courant $T_f=\omega+\ddc g$. On en déduit que
\[\|T_f\|_{X_n}\geq 1-\frac{\deg(f^{\circ n})}{d^n}.\] 
Si $x\in X_n$, alors $g(x)\leq g^{(n)}(x)=-\infty$. Donc la réunion $X=\bigcup X_n$ est contenue dans $\{g=-\infty\}$ et
\[1\geq \|T_f\|_{\{g=-\infty\}}\geq \|T_f\|_X\geq \|T_f\|_{X_n}=1-\frac{\deg(f^{\circ n})}{d^n}.\]
\'Etant donné que $\deg(f^{\circ n})/d^n\to 0$ quand $n\to +\infty$, on a finalement
\[\|T_f\|_{\{g=-\infty\}}=\|T_f\|_X=1.\qedhere\]
\end{proof}

\section{Unicité du courant invariant}

\begin{proposition}\label{prop_unicite}
Soit $f:\P^k\to \P^k$ une application rationnelle de degré $k$ dont le courant de Green dynamique $T_f$ vérifie $\|T_f\|_{\{g=-\infty\}}=1$ avec $g$ un potentiel de $T_f$. Si $T\in \T$ est un courant vérifiant $f^* T = d\cdot T$ et $\|T\|=1$ alors $T=T_f$. 
\end{proposition}

\begin{proof}[Preuve]
Comme nous l'avons rappelé plus haut, $T_f=\omega+\ddc g$ et $T=\omega+\ddc h$ pour des fonctions $\omega$-psh $h$ et $g$ qui vérifient $h\leq g$. Par hypothèse $\|\omega +\ddc g\|_{\{g=-\infty\}}=1$. Le lemme \ref{lemme_psh} ci-dessous implique que $h-g$ est constante et donc que $T=T_f$.
\end{proof}

\begin{lemme}\label{lemme_psh}
Soit $h$ et $g$ deux fonctions $\omega$-psh sur $\P^k$ telles que 
\[h\leq g\quad \text{et}\quad \|\omega+\ddc g\|_{\{g=-\infty\}}=1.\] 
Alors $h=g+c$ avec $c\in \R$. 
\end{lemme}

\begin{proof}[Preuve]
Si $\{g=-\infty\}$ est fermé, la démonstration est immédiate. En effet, en dehors de cet ensemble, la fonction $g$ est lisse et la fonction $h-g$ est psh négative. Elle se prolonge donc en une fonction globalement plurisousharmonique sur $\P^k$, donc constante. 

La difficulté apparait lorsque l'on ne fait pas l'hypothèse que  $\{g=-\infty\}$ est fermé. 
Si l'on travaille localement dans une carte complexe, on peut se donner une fonction lisse $\phi$ telle que $\ddc \phi=\omega$. Les fonctions 
\[u=\phi+g\quad \text{et}\quad v=\phi+h\] sont alors plurisousharmoniques, en particulier sousharmoniques. Si l'on identifie $\C^n$ à $\R^{2n}$, les Laplaciens au sens des distributions de $u$ et $v$ sont des mesures de Riesz $\Delta u$ et $\Delta v$.
Le lieu $\{u\neq -\infty\}$ n'est pas chargé par le courant $\ddc u=\omega+\ddc g$. Il n'est donc pas chargé par la mesure trace de $\ddc u$ qui, à une constante multiplicative près, est égale à $\Delta u$. De plus, $v\leq u$. D'après le lemme \ref{lemme_sousharmonique} ci-dessous, on a donc
\[\Delta v\geq \Delta u.\]
Dans la carte complexe considérée, la fonction $h-g=v-u$ est donc sousharmonique. Comme cela est vrai localement dans n'importe quelle carte complexe, la fonction $h-g$ est globalement plurisousharmonique sur $\P^k$, donc constante. 
%Il s'en suit que globalement sur $\P^k$, on a 
%\[(\omega+\ddc h)\wedge \omega^{k-1}\geq (\omega+\ddc g)\wedge \omega^{k-1}.\]
%\'Etant donné que $\|\omega+\ddc h\|=\|\omega+\ddc g\|$ on en déduit que 
%[(\omega+\ddc h)\wedge \omega^{k-1}= (\omega+\ddc g)\wedge \omega^{k-1}.\]
%Finalement, $\omega+\ddc h=\omega+\ddc g$ et donc $h=g+c$ avec $c\in \R$. 
\end{proof}

\begin{lemme}\label{lemme_sousharmonique}
Soient $m\geq 2$ et soit $B$ la boule unité de $\R^m$. Soient  $u:B\to [-\infty,+\infty)$ et $v:B\to [-\infty,+\infty)$ deux fonctions sousharmoniques et $\Delta u$ et $\Delta v$ les mesures de Riesz sur $B$ associées à $u$ et $v$.  On suppose que 
\[v\leq u\quad \text{et}\quad \Delta u \bigl(\{u\neq -\infty\}\bigr)=0.\] 
Alors $\Delta v \geq \Delta u$. 
\end{lemme}

\begin{proof}[Preuve]
Soit $X\subset B$ un borélien. On doit montrer que $\Delta v(X)\geq \Delta u(X)$. Or, par hypothèse, $\Delta u(X)=\Delta u (Y)$ où $Y=X\cap \{u=-\infty\}$. 
Par propriété de régularité intérieure des mesures boréliennes, on peut trouver une suite de fermés $Y_j\subseteq Y$ tels que $\Delta u(Y_j)\to \Delta u(Y)$ quand $j\to +\infty$. On note $\mu_j={\bf 1}_{Y_j}\cdot \Delta u$ et on choisit pour chaque $j$ un potentiel $u_j:B\to [-\infty,+\infty)$ tel que $\Delta u_j=\mu_j$.\footnote{C'est pour avoir l'existence d'un tel potentiel que l'on travaille dans le cadre des fonctions sousharmoniques et non celui des fonctions plurisousharmoniques.} \'Etant donné que $\Delta u\geq \Delta u_j$, la fonction $u-u_j$ est sousharmonique sur $B$. En particulier, elle est localement majorée. La fonction $u_j$ est harmonique en dehors de $Y_j$. Par conséquent, la fonction $v-u_j$ est sousharmonique sur $B-Y_j$ et localement majorée dans $B$ (par $u-u_j$). Comme $Y_j\subseteq X$ est polaire, il existe une fonction sousharmonique $v_j:B\to [-\infty,+\infty)$ qui prolonge $v-u_j$, c'est-à-dire telle que $v= u_j+v_j$. On a donc 
\[\Delta v = \Delta u_j +\Delta v_j\geq\Delta u_j= \mu_j.\] En particulier
\[\Delta v (X) \geq \Delta v(Y) \geq \mu_j(Y)=\Delta u (Y_j)\underset{j\to +\infty}\longrightarrow \Delta u(Y)=\Delta u(X).\qedhere\]
\end{proof}

\section{Un lemme technique}

Si $X\subseteq \P^k$ est un ensemble mesurable et si $g$ est une fonction $\omega$-psh, on note
\[I_X(g)=\int_X (\omega+\ddc g)\wedge \omega^{k-1}.\]
Dans cette partie, on suppose que $f:\P^k\to \P^k$ et $(f_j:\P^k\to \P^k)_{j\geq 0}$ sont des applications rationnelles de degré $d$ et que $F:\C^{k+1}\to \C^{k+1}$ et $(F_j:\C^{k+1}\to \C^{k+1})$ sont des applications dont les coordonnées sont des polynômes homogènes de degré $d$ telles que
\begin{itemize}
\item $f\circ \pi = \pi \circ F$ et $f_j\circ \pi = \pi \circ F_j$, 
\item $\|F(z)\|<\|z\|^d$ et $\|F_j(z)\|<\|z\|^d$ dès que $\|z\|<1$ et
\item la suite $(F_j)$ converge vers $F$ uniformément sur tout compact de $\C^{k+1}$.
\end{itemize}
On suppose de plus que 
\[F_j(z_j)=z_j\quad \text{avec}\quad 
\C^{k+1}-\{0\}\ni z_j\underset{j\to +\infty}\longrightarrow z\in \C^{k+1}-\{0\}.\]
Comme dans la partie \ref{sec_green}, on définit des suites décroissantes de fonctions $\omega$-psh $(g^{(n)})_{n\geq 0}$ et $(g_j^{(n)})_{n\geq 0}$ qui convergent simplement et dans $L^1(\P^k)$ vers des fonctions $\omega$-psh $g$ et $g_j$. Comme observé plus haut, on a alors
\[g_j\bigl(\pi(z_j)\bigr)=-\log \|z_j\|.\]

Venons-en au lemme clé. 

\begin{lemme}
Si $f$ n'est pas algébriquement stable, alors pour tout  $M\in \R$, on a
\[\lim_{j\to +\infty} I_{\{g_j<-M\}}(g_j) =1.\]
\end{lemme}

\begin{proof}[Preuve]
Fixons $\eps>0$. Soit $n$ suffisamment grand pour que $\deg(f^{\circ n})\leq \eps \cdot d^n$. Alors, comme nous l'avons vu dans la preuve de la proposition \ref{prop_sibony}
\[I_{\{g^{(n)}=-\infty\}}(g)\geq 1-\eps.\]
En particulier pour tout $M'\in \R$, 
\[I_{\{g^{(n)}<-M'\}}(g)\geq 1-\eps.\]

Nous prétendons que la suite de courants $\omega+\ddc g_j$ converge faiblement vers le courant $\omega+\ddc g$. En effet, $g_j\bigl(\pi(z_j)\bigr)=-\log \|z_j\|$ et la suite $(g_j)$ ne peut  pas tendre uniformément vers $-\infty$ sur $\P^k$. 
De plus 
\[g_j\leq g_j^{(n)}\underset{j\to +\infty}\longrightarrow g^{(n)}.\] Par conséquent, toute valeur d'adhérence $h$ de la suite $(g_j)$ dans $L^1(\P^k)$ est majorée par $g^{(n)}$ et comme cela est vrai pour tout $n$, on a $h\leq g$. D'après le lemme \ref{lemme_psh}, on a donc $\omega+\ddc h=\omega+\ddc g$.

La convergence simple de $\omega+\ddc g_j$ vers $\omega+\ddc g$ implique alors que pour $j$ suffisamment grand, 
\[I_{\{g^{(n)}<-M'\} }(g_j) \geq I_{\{g^{(n)}<-M'\} }(g) -\eps \geq 1-2\eps.\]
Finalement, pour $M'>M$, si $j$ est assez grand, alors 
\[\{g^{(n)}<-M'\} \subseteq \{g_j^{(n)}<-M\}\underset{g_j\leq g_j^{(n)}}\subseteq \{g_j<-M\}.\]
Sinon, on pourrait trouver une suite de points $x_j\in \P^k$ avec  $g^{(n)}(x_j)<-M'$ et $g_j^{(n)}(x_j)\geq -M$. Cela signifierait que l'on peut trouver des points $z_j$ de norme 1 dans $\C^{k+1}$ tels que 
$\log \|F^{\circ n}(z_j)\|<-d^n M'$ et $\log \|F_j^{\circ n}(z_j)\|\geq -d^n M$. La suite 
$F_j^{\circ n}$ convergeant uniformément vers $F^{\circ n}$ sur la sphère unité, quitte à extraire une sous-suite de sorte que $z_j$ tende vers $z$, on aurait $-d^n M\leq \log \|F^{\circ n}(z)\|\leq -d^n M'$, d'où une contradiction. Par conséquent, si $j$ est suffisamment grand, 
\[ 1\geq I_{\{g_j<-M\}}(g_j)\geq I_{\{g_j^{(n)}<-M\}}(g_j) \geq I_{\{g^{(n)}<-M'\} }(g_j)\geq 1-2\eps.\qedhere\]
\end{proof}

\section{Existence de courants invariants pluripolaires}

On se place maintenant dans le cadre du théorème. On note $g_t$ la fonction de Green dynamique associée à $f_t$. 
Le lemme précédent implique donc que pour tout $m\geq 1$, l'ensemble
\[O_m = \bigl\{t\in [0,1]~;~I_{\{g_t<-m\}}(g_t)\geq 1-1/m\bigr\}\]
contient un ouvert dense de $[0,1]$ (il contient un voisinage de $\Z$). 
L'intersection $O_m$ est donc un ensemble gras au sens de Baire, de même que 
\[A=\bigcap_{m\geq 0} O_m-\Z.\]

Soit $t\in A$  et notons $f=f_t$ et $g=g_t$. L'application $f$ est algébriquement stable et par construction, pour tout $m\geq 1$, 
\[I_{\{g<-m\}}(g)\geq 1-1/m.\]
En passant à la limite quand $m\to +\infty$ (convergence monotone), on voit que le courant
$T_g$ est pluripolaire~:
\[\|T_g\|_{\{g=-\infty\}} =I_{\{g=-\infty\}}(g)\geq 1=\|T_g\|.\]

\section*{Remerciements}
Je tiens à remercier mes collègues François Berteloot, Dan Coman, Jeffrey Diller, Romain Dujardin, Thomas Gauthier, Vincent Guedj, Frédéric Protin et Ahmed Zeriahi pour les échanges que nous avons eus et qui m'ont permis d'améliorer, voire de corriger des versions préliminaires de cette démonstration.

\end{document}